\documentclass{article}
\usepackage{amsfonts}
\usepackage{amssymb}
\usepackage{amsthm}
\usepackage{latexsym}
\usepackage{graphicx}
\usepackage{amsmath}
\usepackage{epsfig}
\newtheorem{dfn}{Definition}[section]
\newtheorem{ex}{Example}[section]

\newcommand{\bmu}{\boldsymbol{\mu}}

\newcommand{\bx}{\boldsymbol{x}}

\newcommand{\by}{\boldsymbol{y}}

\newcommand{\bzero}{\boldsymbol{0}}
\title{Slicing: Nonsingular Estimation of High Dimensional Covariance Matrices Using Multiway Kronecker Delta Covariance Structures }
\author{Deniz Akdemir \\ Department of Statistics \\ 
  University of Central Florida\\ Orlando, FL 32816 
  }

\begin{document}
\maketitle

\begin{abstract}
Nonsingular estimation of high dimensional covariance matrices is an important step in many statistical procedures like classification, clustering, variable selection an future extraction.  After a review of the essential background material, this paper introduces a technique we call slicing for obtaining a nonsingular covariance matrix of high dimensional data. Slicing is essentially assuming that the data has Kronecker delta covariance structure. Finally, we discuss the implications of the results in this paper and provide an example of classification for high dimensional gene expression data.\end{abstract}

%\keywords{Normal Distribution, Multivariate Distribution, Matrix Variate Normal Distribution, Array Variate Random Variable, Array Variate Normal Distribution, Classification, Dimension Reduction, High Dimensional Data}
%\ams{Primary 62H10, Secondary 62H05.}

\maketitle

\section{Intoduction}

The advances in data collection methods and the increase in data storage and processing capabilities has led to data sets that are not suitable for analysis with the classical statistical approaches. For example, through DNA micro array techniques, the expression levels of millions of genes can easily be obtained.  However, usually, the number of observations (the sample size) is much less than the number of expression levels observed. This is the characteristic of many recent data sets in bioinformatics, signal processing, and many other fields of science. The number of variables (p) is much higher than the number of observations (N) (i.e., $N<<p$).

It is well known that when $N<p$ the usual sample covariance matrix will be singular. Many methods in statistics, like clustering and classification depends on estimating the inverse of the covariance matrix. For small samples and especially when $N<<p$  This becomes a major problem when the we need to obtain a the inverse of the covariance matrix.  

The technique slicing, which we will discuss in detail in this paper, is essentially obtaining estimates of the covariance matrix under the assumption that assuming that the p-dimensional observations are realizations from a multivariate distribution with a certain Kronecker delta structure. Slicing is appropriate when the number of observations in the sample is much less than the number of variables.because by choosing a Kronecker structure for the covariance a great deal of decrease in the number of parameters is obtained. By using $2$-way, $3$-way, and in general $i$-way Kronecker structures for the covariance matrix, we can obtain nonsingular estimates of the covariance matrix when $N<<p.$

While developing slicing, we have used the concept of array variate normal variable with multiway Kronecker delta structure obtained by using the rules of multi linear algebra. In Section 2, we will first review array algebra as its discussed in \cite{rauhala1974array},  \cite{rauhala1980introduction}, Blaha \cite{blaha1977few}. The array variate normal model with Kronecker delta structure and estimation of its parameters are also discussed in Section 2. In Section 3, we describe slicing in detail, provide the results from various simulations and apply the technique to high dimensional gene expression data.  

\section{Array Algebra and Array Variate Normal Random Variable}
\subsection{Array Algebra}

In this paper we will only study arrays with real elements. We will write $\widetilde{X}$ to say that $\widetilde{X}$ is an array. When it is necessary we can write the dimensions of the array as subindices, e.g., if $\widetilde{X}$ is a $m_1 \times m_2\times m_3 \times m_4$ dimensional array in $R^{m_1\times m_2 \times \ldots \times m_i}$, then we can write $\widetilde{X}_{m_1 \times m_2\times m_3 \times m_4}.$  To refer to an element of an array $\widetilde{X}_{m_1 \times m_2\times m_3 \times m_4},$ we write the position of the element as a subindex to the array name in parenthesis, $(\widetilde{X})_{r_1r_2r_3r_4}.$ 
 
We will now review some basic principles and techniques of multi linear algebra. These results and their proofs can be found in Rauhala \cite{rauhala1974array},  \cite{rauhala1980introduction}  and \cite{blaha1977few}. 

\begin{dfn}\emph{Inverse Kronecker product} of two matrices $A$ and $B$ of dimensions $p\times q$ and $r \times s$ correspondingly is written as $A\otimes^i B$ and is defined as $A\otimes^i B=[A(B)_{jk}]_{pr\times qs}=B\otimes A,$ where $'\otimes'$ represents the ordinary Kronecker product.\end{dfn}

The following properties of the inverse Kronecker product are useful:\begin{itemize}\item $\bzero \otimes^i A= A \otimes^i \bzero=\bzero.$ \item $(A_1+A_2)\otimes^i B=A_1\otimes^i B+ A_2 \otimes^i B.$ \item $A \otimes^i (B_1+ B_2)=A \otimes^i B_1+ A \otimes^i B_2.$ \item  $\alpha A \otimes^i \beta B= \alpha \beta A \otimes^i B.$ \item $(A_1 \otimes^i B_1)(A_2 \otimes^i B_2)= A_1A_2 \otimes^i B_1B_2.$ \item  $(A \otimes^i B)^{-1}=(A^{-1} \otimes^i B^{-1}).$ \item  $(A \otimes^i B)^{+}=(A^{+} \otimes^i B^{+}),$ where $A^{+}$ is the Moore-Penrose inverse of $A.$ \item $(A \otimes^i B)^{-}=(A^{-} \otimes^i B^{-}),$ where $A^{-}$ is the $l$-inverse of $A$ defined as $A^{-}=(A'A)^{-1}A'.$\item If $\{\lambda_i\}$ and $\{\mu_j\}$ are the eigenvalues with the corresponding eigenvectors $\{\bx_i\}$ and $\{\by_j\}$ for matrices $A$ and $B$ respectively, then $A\otimes^i B$ has eigenvalues $\{\lambda_i\mu_j\}$ with corresponding eigenvectors $\{\bx_i\otimes^i\by_j\}.$\item Given two matrices $A_{n\times n}$ and $B_{m\times m}$ $|A\otimes^i B|=|A|^m|B|^n,$ $tr(A\otimes^i B)=tr(A)tr(B).$\item $A\otimes^i B=B\otimes A=U_1 A\otimes B U_2,$ for some permutation matrices $U_1$ and $U_2.$ 
\end{itemize}

It is well known that a matrix equation $$AXB'=C$$ can be rewritten in its mono linear form as \begin{equation}\label{eqmnf}A\otimes^i B vec(X)=vec(C).\end{equation} Furthermore, the matrix equality $$A\otimes^i B XC'=E$$ obtained by stacking equations of the form (\ref{eqmnf}) can be written in its mono linear form as $$(A\otimes^i B \otimes^i C) vec(X)=vec(E).$$ This process of stacking equations could be continued and R-matrix multiplication operation introduced by Rauhala \cite{rauhala1974array} provides a compact way of representing these equations in array form:

\begin{dfn}\emph{R-Matrix Multiplication} is defined element wise: 

 $$((A_1)^1 (A_2)^2 \ldots (A_i)^i\widetilde{X}_{m_1 \times m_2 \times \ldots \times m_i})_{q_1q_2\ldots q_i}$$ $$=\sum_{r_1=1}^{m_1}(A_1)_{q_1r_1}\sum_{r_2=1}^{m_2}(A_2)_{q_2r_2}\sum_{r_3=1}^{m_3}(A_3)_{q_3r_3}\ldots \sum_{r_i=1}^{m_i}(A_i)_{q_ir_i}(\widetilde{X})_{r_1r_2\ldots r_i}.$$ \end{dfn}

R-Matrix multiplication generalizes the matrix multiplication (array multiplication in two dimensions)to the case of $k$-dimensional arrays. The following useful properties of the R-Matrix multiplication are reviewed by Blaha \cite{blaha1977few}:
\begin{itemize}
\item $(A)^1B=AB.$ 
\item $(A_1)^1(A_2)^2C=A_1CA'_2.$ 
\item $\widetilde{Y}=(I)^1(I)^2\ldots (I)^i \widetilde{Y}.$ 
\item \small $((A_1)^1 (A_2)^2\ldots (A_i)^i)((B_1)^1(B_2)^2\ldots (B_i)^i)\widetilde{Y}= (A_1B_1)^1(A_2B_2)^2\ldots(A_iB_i)^i\widetilde{Y}.$ \normalsize
\end{itemize}

The operator $rvec$ describes the relationship between $\widetilde{X}_{m_1 \times m_2 \times \ldots m_i}$ and its mono linear form $\bx_{m_1m_2\ldots m_i\times 1}.$ 
\begin{dfn}\label{def:rvec} $rvec( \widetilde{X}_{m_1 \times m_2 \times \ldots m_i})=\bx_{m_1m_2\ldots m_i\times 1}$ where $\bx$ is the column vector obtained by stacking the elements of the array $\widetilde{X}$ in the order of its dimensions; i.e., $(\widetilde{X})_{j_1 j_2 \ldots j_i}=(\bx)_j$ where $j=(j_i-1)n_{i-1}n_{i-2}\ldots n_1+(j_i-2)n_{i-2}n_{i-3}\ldots n_1+\ldots+(j_2-1)n_1+j_1.$\end{dfn}

Let $\widetilde{L}_{m_1 \times m_2 \times\ldots m_i}=(A_1)^1(A_2)^2\ldots(A_i)^i\widetilde{X}$ where $(A_j)^j$ is an $m_j\times n_j$ matrix for $j=1,2,\ldots,i$ and $\widetilde{X}$ is an $n_1\times n_2\times\ldots\times n_i$ array. Write $\mathbf{l}=rvec(\widetilde{L})$ and $\bx=rvec(\widetilde{X}).$ Then, $\mathbf{l}=A_1\otimes^iA_2\otimes^i\ldots\otimes^i A_i\bx.$ Therefore, there is an equivalent expression of the array equation in mono linear form.

\begin{dfn}{} The square norm of $\widetilde{X}_{m_1 \times m_2 \times\ldots m_i}$ is defined as $$\|\widetilde{X}\|^2=\sum_{j_1=1}^{m_1}\sum_{j_2=1}^{m_2}\ldots\sum_{j_i=1}^{m_i}((\widetilde{X})_{j_1j_2\ldots j_i})^2.$$ \end{dfn}

\begin{dfn}{} The distance of $\widetilde{X_1}_{m_1 \times m_2 \times\ldots m_i}$ from $\widetilde{X_2}_{m_1 \times m_2 \times\ldots m_i}$ is defined as $$\sqrt{\|\widetilde{X_1}-\widetilde{X_2}\|^2}.$$ \end{dfn}

\begin{ex} Let $\widetilde{Y}=(A_1)^1 (A_2)^2\ldots (A_i)^i\widetilde{X}+\widetilde{E}.$ Then $\|\widetilde{E}\|^2$ is minimized for $\widehat{\widetilde{X}}=(A_1^{-})^1(A_2^{-})^2\ldots(A_i^{-})^i\widetilde{Y}.$ \end{ex}

\subsection{Array Variate Normal Distribution}

\begin{dfn}(\cite{DenizGuptaJAS})\label{modkroncov} Let $A_1, A_2,\ldots,A_i$ are non singular matrices of orders $m_1, m_2,\ldots, m_i$ and let $\widetilde{M}$ be an $m_1 \times$ $m_2$ $\times \ldots$ $\times m_i$  dimensional constant array.  Then the pdf of  array normal random variable $\widetilde{X}$ with Kronecker delta covariance structure is given by \small \begin{equation}\label{eq:densityarn}\phi(\widetilde{X}; \widetilde{M},A_1,A_2,\ldots A_i)=\frac{\exp{(-\frac{1}{2}\|{(A_1^{-1})^1 (A_2^{-1})^2 \ldots (A_i^{-1})^i(\widetilde{X}-\widetilde{M})}\|^2)}}{(2\pi)^{m_1m_2\ldots m_i/2}|A_1|^{\prod_{j\neq 1}{m_j}} |A_2|^{\prod_{j\neq 2}{m_j}} \ldots |A_i|^{\prod_{j\neq i}{m_j}}}.\end{equation} \normalsize   \end{dfn}

Distributional properties of a array normal variable with density in the form of Theorem \ref{modkroncov} can obtained by using the equivalent  mono linear representation. The moments, the marginal and conditional distributions, independence of variates should be studied considering the equivalent mono linear form of the array variable and the well known properties of the multivariate normal random variable. 

\begin{dfn} For the $m_1 \times m_2 \times\ldots \times m_i$ dimensional array variate random variable $\widetilde{X},$ the principal components are defined as the principal components of the  $d=m_1m_2\ldots m_i$-dimensional random vector $rvec(\widetilde{X}).$ \end{dfn} 

The main statistical problem is the estimation of the covariance of $rvec(\widetilde{X}),$ its eigenvectors and eigenvalues for small sample sizes. 

\subsection{Estimation}

In this section we provide an heuristic method of estimating the model parameters. The optimality of these estimators are not proven but merely checked by simulation studies. Inference about the parameters of the model in Theorem \ref{modkroncov}  for the matrix variate case has been considered in the statistical literature (\cite{roy2003tests}, \cite{roy2008likelihood}, \cite{lu2005likelihood},\cite{srivastava2008models}, etc...). In these papers, the unique maximum likelihood estimators of the parameters of the model in Theorem \ref{modkroncov} for the matrix variate case are obtained under different assumptions for the covariance parameters. Some classification rules based on the matrix variate observations with Kronecker delta covariance structures have been studied in \cite{roy2009classification}, and also in \cite{krzysko2009discriminant}.

The model in Theorem \ref{modkroncov} the way it is stated is  unidentifiable. However, this problem can easily be resolved by putting restrictions on the covariance parameters. The approach we take is to assume that $j-1$ of the last diagonal elements of matrices $A_jA'_j$  are equal to $1$ for $j=1,2,\ldots,i.$ The Flip-Flop Algorithm is proven to attain the maximum likelihood estimators of the parameters of two dimensional array variate normal distribution  \cite{srivastava2008models}. 

The following is similar to the flip flop algorithm. First, assume $\{\widetilde{X}_1,$ $\widetilde{X}_2,$ $\ldots,$ $\widetilde{X}_N\}$ is a random sample from a $N(\widetilde{M},A_1,A_2,\ldots A_i)$ distribution with $j-1$ of the last diagonal elements of matrices $A_jA'_j$ equal to $1$ for $j=1,2,\ldots,i.$ Further, we assume that all $A_j'$s are square positive definite matrices of rank at least $j.$ Finally, assume that we have $N\prod_{j=1}^{i}m_j>m_r^2$ for all $r=1,2, \ldots,i.$

Algorithm for estimation:
\begin{enumerate}
	\item Estimate $\widetilde{M}$ by $\widehat{\widetilde{M}}=\frac{1}{N}\sum_{l=1}^N \widetilde{X}_l,$ and obtain the centered array observations $\widetilde{X}_l^c=\widetilde{X}_l-\widehat{\widetilde{M}}$ for $l=1,2,\ldots, N.$ 
	\item Start with initial estimates of $A_2, A_3, \ldots, A_i.$ 
	\item On the basis of the estimates of $A_2, A_3, \ldots, A_i$ calculate an estimate of $A_1$ by first scaling the array observations using $$\widetilde{Z}_l=(I)^1(A_2^{-1})^2, (A_3^{-1})^3, \ldots, (A_i^{-1})^i\widetilde{X}_l^c,$$ and then calculating the square root of covariance along the $1$st dimension of the arrays $\widetilde{Z}_l,$ $l=1,2,\ldots, N.$   
	\item On the basis of the most recent estimates of the model parameters, estimate $A_j$ $j=2,\ldots,i.$ by first scaling the array observations using $$\widetilde{Z}_l=(A_1^{-1})^1(A_2^{-1})^2, \ldots (A_{j-1}^{-1})^{j-1} I (A_{j+1}^{-1})^{j+1}\ldots(A_i^{-1})^i\widetilde{X}_l^c,$$ and then calculating the square root of covariance along the jth dimension of the arrays $\widetilde{Z}_l$'s for $j=2,\ldots, i.$	Scale  the estimate of $A_jA'_j$ so that the last $j-1$ diagonal elements are equal to $1.$ 
	\item Repeat steps 3 and 4 until convergence is attained. 
\end{enumerate}

Let  $\widetilde{X}_l,$ $l=1,2,..,N$ be a random sample for the array variate random variable $\widetilde{X}.$ Let $p=m_1 m_2 \ldots m_i.$ When $N<p,$ it is well known that the usual covariance estimator for $rvec(\widetilde{X})$ will be singular with probability one. Therefore, when $N<p,$ there is no consistent estimator of the covariance of $rvec(\widetilde{X})$ under the unstructured covariance assumption.

On the other hand, if we assume that the covariance matrix has Kronecker delta structure, we can obtain a nonsingular estimate of the covariance structure with the methods developed in this section. The condition on the sample size is relaxed considerably. If we have $pN>m_r^2$ for $r=1,2, \ldots,i$ and the assumptions stated before the algorithm for estimation of the parameters of this model hold, then the estimator of the covariance matrix is nonsingular. When the covariance does not have Kronecker structure, the estimate obtained here could be used as regularized nonsingular estimate of the covariance.

\begin{ex}\label{exsim1}Let \\  \begin{center} \small  $A_1= \left( \begin{array}{cc}
4 & 1  \\
1 & 2 \end{array} \right)^{1/2}, 
$ $A_2= \left( \begin{array}{ccc}
3 & 0 & -1 \\
0 & 2 & 0 \\
-1 & 0 & 1 \end{array} \right)^{1/2},$ \normalsize  \\
\end{center} 
and 
\begin{center} \small $A_3=\left( \begin{array}{ccc}
4 & 0 & 1 \\
0 & 1 & 0 \\
1 & 0 & 1 \end{array} \right)^{1/2}. 
$  \end{center} \normalsize Also, let $\widetilde{M}$ be the $0$ array of dimensions $2 \times 3 \times 3.$ The following are the estimates of $A_1$, $A_2$ and $A_3$ based on a random sample of size $100$ from the $N(A_1, A_2,A_3, \widetilde{M}).$ \\ \small \begin{center}$\widehat{A}_1= \left( \begin{array}{cc}
2.76 & 0.48  \\
0.48 & 1.13 \end{array} \right)^{1/2}, 
$ $\widehat{A}_2= \left( \begin{array}{ccc}
4.69 & 0.25 & -0.52 \\
0.25 & 2.70 & -0.04 \\
-0.52 & -0.04 & 1 \end{array} \right)^{1/2},$ \\ 
\normalsize \end{center} and \\ \small \begin{center} $\widehat{A}_3=\left( \begin{array}{ccc}
4.27 & -0.13 & 0.55 \\
-0.13 & 1 & 0.02 \\
0.55 & 0.02 & 1 \end{array} \right)^{1/2}.$ \normalsize  \\
\end{center}
The left plot in Figure \ref{figuresim} compares the estimated eigenvalues to the true eigenvalues for this example.
\end{ex}

\section{Slicing}

A vector $\bx$ of dimension $p$ can be sliced into $p/m_1=m_2$ pieces and organized into a matrix of order $p=m_1\times m_2$ for some natural numbers $m_1$ and $m_2.$ Or, in general, the same vector can be organized in an array of dimension $p=m_1\times m_2\times \ldots \times m_i$  for some natural numbers $m_1,$ $m_2, \ldots,$ $m_i.$ Once we slice the data and reorganize it in array form, we can pretend that this array data was generated from the model in Theorem \ref{modkroncov}. We require that the additional assumptions stated before the algorithm for estimation of the parameters of this model hold. A nonsingular estimate of the covariance matrix $\Lambda$ of the $p-$dimensional vector variate random variable can be obtained by using the estimators from this algorithm and using $\hat{\Lambda}=(\hat{A}_1\otimes^i \hat{A}_2\otimes^i\hat{A}_i)(\hat{A}_1\otimes^i \hat{A}_2\otimes^i\hat{A}_i)'$ 

That we do not have to assume any covariance components are zero is the main difference and advantage of this regularization method to the usual shrinkage methods like lasso \cite{friedman2008sparse}.

If $\{\lambda(A_r)_{r_j}\}$ are the $m_j$ eigenvalues of $A_rA'_r$ with the corresponding eigenvectors $\{(\bx_r)_{r_j}\}$ for $r=1,2,\ldots,i$ and $r_j=1,2,\ldots,m_r,$ then $(A_1\otimes^i A_2\otimes^iA_i)(A_1\otimes^i A_2\otimes^iA_i)'$ will have  eigenvalues $\{\lambda(A_1)_{r_1}\lambda(A_2)_{r_2}\ldots\lambda(A_i)_{r_i}\}$ with corresponding eigenvectors $\{(\bx_1)_{r_1}\otimes^i(\bx_2)_{r_2}\otimes^i \ldots \otimes^i(\bx_i)_{r_i} \}.$  By replacing  $A_r$ by their estimators, we estimate the eigenvalues and eigenvectors of the covariance of $rvec(\widetilde{X})$ using this relationship. Since each eigenvector is a Kronecker product of smaller components the reduction in dimension obtained by this approach is larger than the one that could be obtained using the ordinary principle components on the ordinary sample covariance matrix.

\begin{ex}Let $\bx\sim N_{12}(\bmu=0, \Lambda).$. We illustrate slicing for $i=2,$ $m_1=3$ and $m_2=4.$ $N=5,10, 15,20,...,45,50$ sets of $N$ observations were generated and $\Lambda$ was estimated using $\widehat{\Lambda}=\widehat{A}_1\otimes^i \widehat{A}_2$ assuming the model in Theorem \ref{modkroncov}. We repeated the whole experiment $5$ times. The results are summarized in Figure \ref{fig:slicingN}. The covariance matrix in the left figure is the identity matrix. In the center figure we have $\Lambda$ that has the same order Kronecker delta covariance structure as the slicing, the components of $\Lambda$ are unstructured and generated randomly. The right figure is the case where $\Lambda$ is a randomly generated unstructured covariance matrix. Slicing has as a regularization effect that shrinks the eigenvalues towards each other.

\begin{figure}[htbp]
	\centering
		\includegraphics[width=1\textwidth]{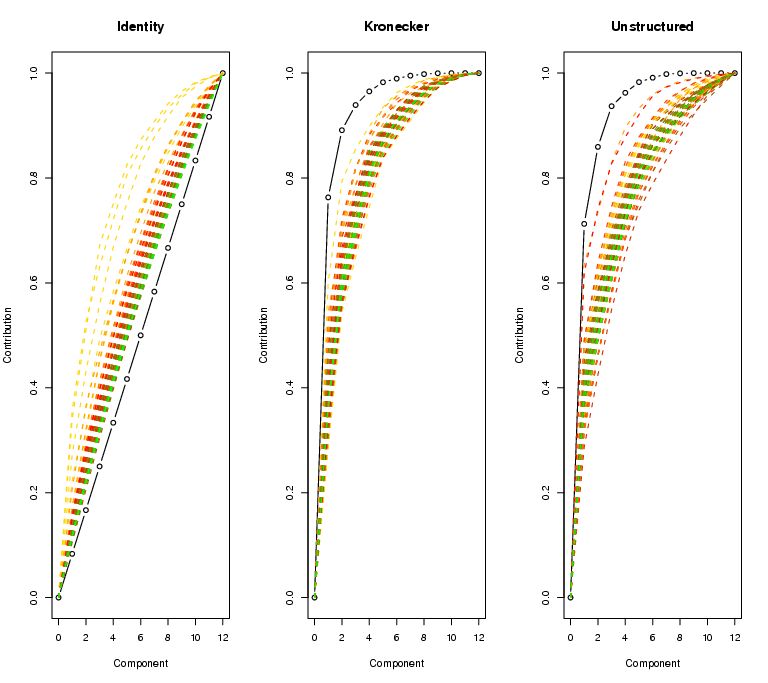}
	\caption{The true covariance structure is represented with the black line. The estimated covariance structure is denoted by different colors according to the sample size. Yellow colors are for small sample sizes, red for moderate sample sizes and green for the larger samples.  The covariance matrix in the left figure is the identity matrix. As N gets larger the estimators of the eigenvalues approach the true values. In this case the estimator seem to be consistent. In the center figure $\Lambda$ has the same order Kronecker delta covariance structure as the slicing, each of the components of $\Lambda$ are unstructured and generated randomly. The right figure is the case where $\Lambda$ is a randomly generated unstructured covariance matrix. In these last two cases the estimator has a bias that does not decrease on the average with increasing sample sizes, however the variance of the estimator decreases as the sample size increases.}
	\label{fig:slicingN}
\end{figure}

\end{ex}

\begin{ex}\label{alon1} The Alon colon data set \cite{alon1999broad} have expression measurements on 2000 genes and $n_1=40$ tumor tissues and $n_2=22$ normal tissue samples. We will compare the means of the normal and tumor tissue samples. We assume first that normal and tumor tissues have the same covariance $\Lambda,$ a $2000\times 2000$ positive definite matrix. We slice each of the $n=62$ observations into a $40\times 50$ matrix and estimate $\Lambda$ with $\widehat{\Lambda}=\widehat{A}_1\otimes^i \widehat{A}_2$ assuming the model in Theorem \ref{modkroncov} holds. For testing the equality of the means, we calculate the $F^+$ statistic proposed in \cite{kubokawa2008estimation} replacing their estimator of the inverse of covariance matrix $\Lambda$ with the inverse of $\widehat{\Lambda}$: \[F^{+}=\frac{2000-(62-1)+1}{(62-1)^2}\left(\frac{1}{40}+\frac{1}{22}\right)^{-1}(\bar{\bx}_1-\bar{\bx}_2)'\widehat{\Lambda}^{-1}(\bar{\bx}_1-\bar{\bx}_2)=3023.273.\]  Using the sampling distribution $F_{r,n-r}$ proposed in \cite{kubokawa2008estimation} assuming that the rank $r$ of $\widehat{\Lambda}$ is $62-1,$ the p-value is calculated as $0.01445.$ Thus, the hypothesis of equality of the means is rejected. 

\end{ex}

\begin{figure}[htbp]
  \begin{center}
    \includegraphics[width=.4\textwidth]{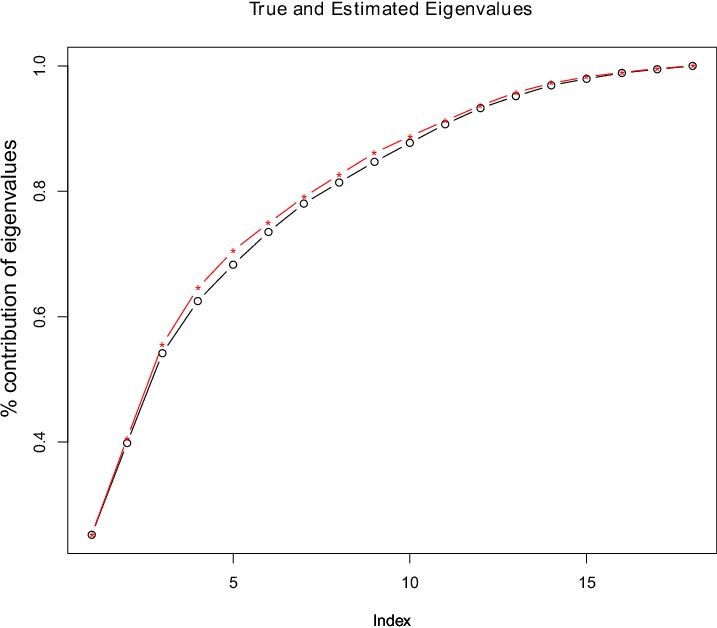}
    \includegraphics[width=.4\textwidth]{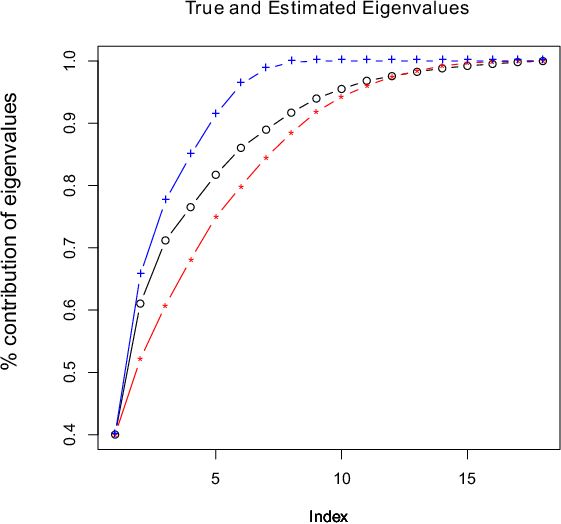}
  \end{center}
  \caption{The left plot in Figure \ref{figuresim} compares the estimated eigenvalues to the true eigenvalues (Example \ref{exsim1}). The right plot in Figure \ref{figuresim} compares the estimated eigenvalues obtained under different assumptions to the true eigenvalues (Example \ref{exsim2}). The red $*$'s represent the true values,  black $\circ$'s are for estimates under Kronecker delta covariance assumption and the blue $+$'s are for estimates under unrestricted covariance assumption.}
  \label{figuresim}
\end{figure}

\begin{ex}\label{exsim2} $N=10$   i.i.d. observations  from a \small $N_{12}(\bmu=0, \Lambda),$ \normalsize  \normalsize distribution are generated for a randomly generated unstructured nonsingular covariance matrix  $\Lambda$. The right plot in Figure \ref{figuresim} compares the estimated eigenvalues obtained by slicing this data into a $2 \times 3 \times 3$ array with the ordinary sample covariance.\end{ex}

\begin{ex} In this example, we will use the heatmap of the true and estimated covariance matrices under different scenarios to see that slicing gives a reasonable description of the variable variances and covariances. In Figure \ref{fig:heatmap120var15t8identity} the true covariance matrix is a $120\times 120$ identity matrix, we estimate this covariance matrix for $N=10, 50,$ and $100$ independent sets of random samples by using $15\times 8$ slicing. In Figure \ref{fig:heatmap120var15t8kron} the true covariance matrix is a $120\times 120$ block diagonal matrix with Kronecker delta structure. Finally, in Figure \ref{fig:heatmap120var15t8unstruct} the true covariance is a matrix with 4 way Kronecker structure. Convergence of the estimators is observed even when $p>>N.$ 

\begin{figure}[htbp]
	\centering
		\includegraphics[width=.5\textwidth]{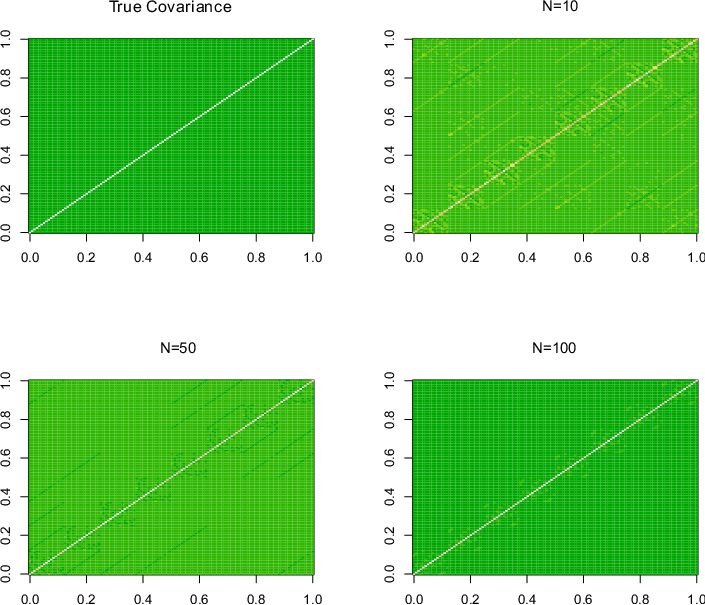}
	\caption{The true covariance matrix is a $120\times 120$ identity matrix with Kronecker delta structure. We estimate this covariance matrix for $N=10, 50,$ and $100$ independent sets of random samples by using $15\times 8$ slicing.}
	\label{fig:heatmap120var15t8identity}
\end{figure}

\begin{figure}[htbp]
	\centering
		\includegraphics[width=.5\textwidth]{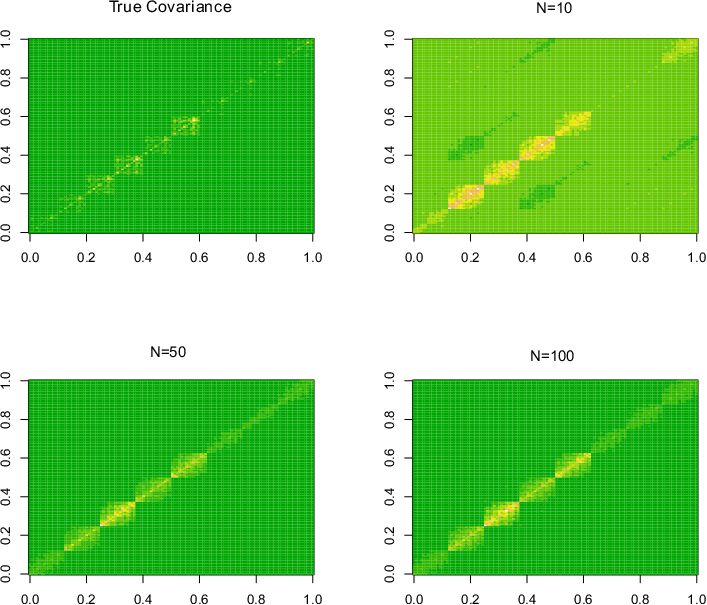}
	\caption{The true covariance matrix is a $120\times 120$ block diagonal matrix with Kronecker delta structure. We estimate this covariance matrix for $N=10, 50,$ and $100$ independent sets of random samples by using $15\times 8$ slicing.}
	\label{fig:heatmap120var15t8kron}
\end{figure}

\begin{figure}[htbp]
	\centering
		\includegraphics[width=.5\textwidth]{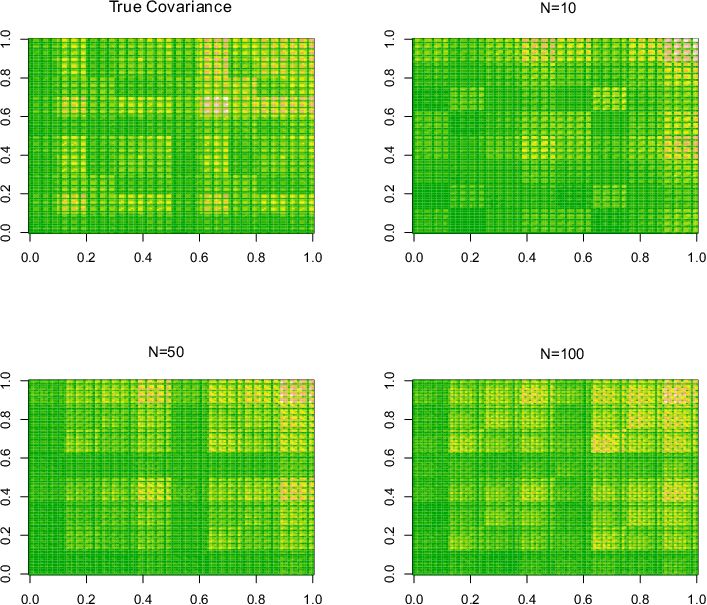}
	\caption{The true covariance matrix is $120\times 120$ 4-way Kronecker delta structured matrix. We estimate this covariance matrix for $N=10, 50,$ and $100$ independent sets of random samples by using $15\times 8$ slicing.}
	\label{fig:heatmap120var15t8unstruct}
\end{figure}

\end{ex}

\begin{ex} 
\begin{figure}[t]
	\centering
		\includegraphics[width=.5\textwidth]{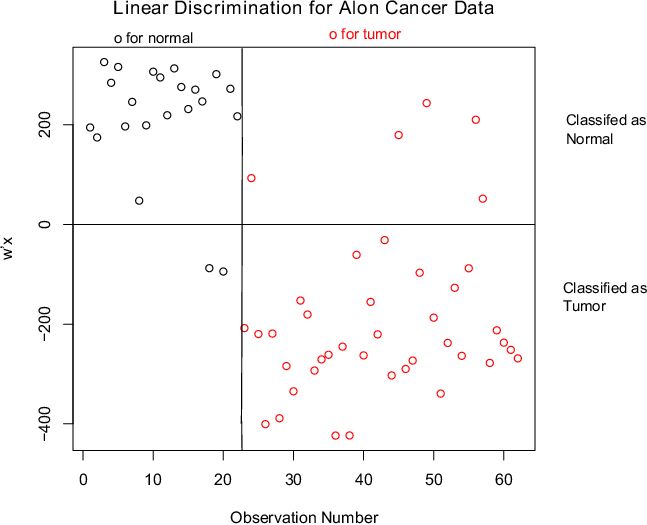}
    \caption{Linear discriminant analysis for the Alon colon data set \cite{alon1999broad}. An observation $\bx$ was classified as ''normal'' if $\bx'\mathbf{w}>0,$ otherwise as ''tumor''. Misclassification rate is $\%11.3.$}
		\label{fig:colonclassification}
\end{figure}

We have used the Fisher's linear discriminant analysis for the Alon colon data set \cite{alon1999broad}. The linear discriminant function was calculated using $\mathbf{w}=\widehat{\Lambda}^{-1}(\bar{\bx}_1-\bar{\bx}_2)$ where $\widehat{\Lambda}$ is the covariance estimate from Example \ref{alon1}. An observation $\bx$ was classified as ''normal'' if $\bx'\mathbf{w}>0,$ otherwise as ''tumor''. Figure \ref{fig:colonclassification} summarizes our findings. Misclassification rate is $\%11.3.$ \end{ex}

In practice, how slicing is done matters. For example,  a $24$ dimensional vector could be sliced as $2\times 12,$ $3 \times 8,$ $4 \times 6,$ or  $2\times 3 \times 4, $ etc. In addition, the permutation of the variables will effect the estimators. As was discussed earlier slicing obtains dimension reduction by writing the covariance matrix into separable components and we perceive that more parsimonious models can be obtained by, for example, proposing a reduced rank mean for the array variable obtained after slicing. Yet another direction would be estimating each component of the covariance structure sparsely by using  a penalty approach like the one used in  \cite{friedman2008sparse}. These issues and improvements are important and will be dealt with in detail in a different article. In the following,we will use the GLASSO \cite{kim2005glasso} package which implements the shrinkage estimator of covariance \cite{friedman2008sparse} matrices will be used in conjunction with the flip flop algorithm. In practice, each of the components of the covariance structure could be penalized to individually to obtain very sparse nonsingular covariance estimates. This is important for variable selection.

\begin{ex}(Sparse Slicing with GLASSO:)

\begin{figure}[htbp]
	\centering
		\includegraphics[width=1\textwidth]{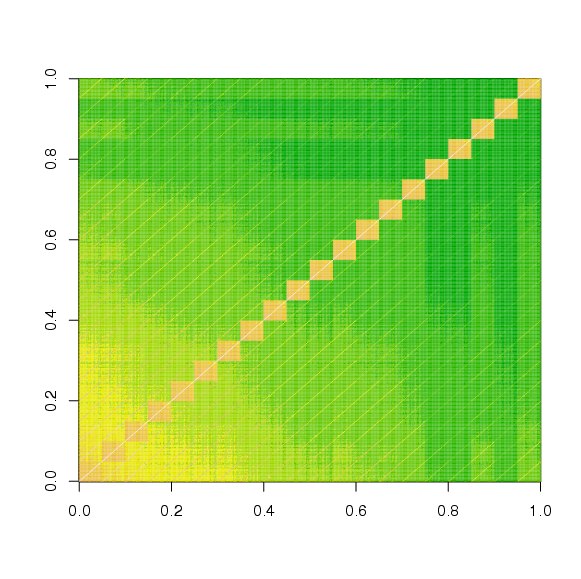}
	\caption{Slicing with GLASSO: The expression levels in this dataset were ordered with respect to their variances. For the samples of expression levels from normal and tumor tissues, high correlation values (lighter colors in the heatmap) are only observed for the expression levels that have high variance. Low variance components have little correlation among each other but they might be mildly correlated with the high variance expression levels.}
	\label{fig:lastcolon}
\end{figure}

We insert the GLASSO algorithm of \cite{friedman2008sparse} at the 4th step of the estimation algorithm from Section 2.3, just before scaling of the matrix. The heatmap of the estimated correlation matrix for the first 500 components of the Alon colon data set obtained by using two way slicing ($20\times 25$) and applying GLASSO to the components are given in Figure \ref{fig:lastcolon}. The shrinkage parameters for GLASSO should selected by the aid of a model selection technique. Here, the values of these parameters are identified tentatively. The expression levels in this dataset were ordered with respect to their variances. For the samples of expression levels from normal and tumor tissues, high correlation values (lighter colors in the heatmap) are only observed for the expression levels that have high variance. Low variance components have little correlation among each other but they might be mildly correlated with the high variance expression levels. The linear discrimination of the groups on the 500 high variance expression levels result in $12.9\%$  false classification rate.

\end{ex}

\bibliographystyle{plain}
\bibliography{arrayref}
\end{document}